\documentclass[11pt,oneside,notitlepage]{amsart}
\usepackage[french]{babel}
\usepackage[latin1]{inputenc}
\usepackage{mathpazo}
\usepackage[symbol,perpage]{footmisc}
\usepackage{xspace}

\vfuzz \hfuzz
\let\Person\bsc

\let\Title\textit

\begin{document}

\begin{center}
\textsc{\Large\bf Sur la convergence des séries trigonométriques qui servent à représenter une fonction arbitraire entre des limites
données.}\footnote{Originally published in Journal f\"ur die reine und
angewandte Mathematik Vol. 4 (1829) p. 157--169.
Transcribed by Ralf Stephan,
eMail: \texttt{mailto:ralf@ark.in-berlin.de}. The scanned images are available
at 
\texttt{http://gdz.sub.uni-goettingen.de/no\_cache/dms/load/img/?IDDOC=270602}}\\[5pt]
(Par Mr. \textit{Lejeune-Dirichlet}, prof. de mathém.)\\[15pt]
\end{center}

\bigskip

MSC-Class: 51-03 01A16

\bigskip

Les séries de sinus et de cosinus, au moyen desquelles on peut représenter
une fonction arbitraire dans un intervalle donné, jouissent entre
autres propriétés remarquables aussi de celle d'être convergentes. Cette
propriété n'avait pas échappée au géomètre illustre qui a ouvert une
nouvelle carrière aux applications de l'analyse, en y introduisant la manière
d'exprimer les fonctions arbitraires dont il est question; elle se
trouve énoncée dans le Mémoire qui contient ses premières recherches
sur la chaleur. Mais personne, que je sache, n'en a donné jusqu'a présent
une démonstration générale. Je ne connais sur cet objet qu'un travail
dû à M.~\Person{Cauchy} et qui fait partie des Mémoires de l'Académie
des sciences de Paris pour 1'année 1823. L'auteur de ce travail avoue
lui même que sa démonstration se trouve en défaut pour certaines fonctions
pour lesquelles la convergence est pourtant incontestable. Un examen
attentif du Mémoire cité m'a porté à croire que la démonstration
qui y est exposée n'est pas même suffisante pour les cas auxquels l'auteur
la croit applicable. Je vais, avant d'entrer en matière, énoncer en
peu de mots les objections auxquelles la démonstration de M.~\Person{Cauchy} me
parait sujette. La marche que ce géomètre célèbre suit dans cette recherche
exige que l'on considère les valeurs que la fonction $\phi(x)$
qu'il s'agit de développer,
obtient, lorsqu'on y remplace la variable $x$ par une quantité de
la forme $u+v\sqrt{-1}$. La considération de ces valeurs semble étrangère
à la question et l'on ne voit d'ailleurs pas bien ce que l'on doit entendre
par le résultat d'une pareille substitution lorsque la fonction dans laquelle
elle a lieu, ne peut pas être exprimée par une formule analytique. Je
présente cette objection avec d'autant plus de confiance, que l'auteur me
semble partager mon opinion sur ce point. Il insiste en effet dans plusieurs
de ces ouvrages sur la nécessité de définir d'une manière précise
le sens que l'on attache à une pareille substitution même lorsqu'elle est
faite dans une fonction d'une loi analytique régulière; on trouve surtout
dans le Mémoire qu'il a inseré dans le 19\up{ième} cahier du journal
polytechnique pag.~567 et suiv.\ des remarques sur les difficultés que font 
naître les quantités imaginaires placées sous des signes de fonctions
arbitraires.
Quoi qu'il en soit de cette première observation, la démonstration de
M.~\Person{Cauchy} donne encore lieu à une autre objection qui parait ne 
laisser aucun doute sur son insuffisance. La considération des quantités
imaginaires conduit l'auteur à un résultat sur le décroissement des termes
de la série, qui est loin de prouver que ces termes forment une
suite convergente. Le résultat dont il s'agit peut être énoncé comme il
suit, en supposant que l'intervalle que l'on considère, s'etende depuis zéro
jusqu'à $2\pi$.
\begin{quote}
\og Le rapport du terme dont le rang est $n$, à la quantité
$A\frac{\sin nx}{n}$
($A$ désignant une constante determinée dépendante des valeurs extrêmes
de la fonction) diffère de l'unité prise positivement d'une quantité qui 
diminue indéfiniment, à mesure que $n$ devient plus grand.\fg
\end{quote}

De ce résultat et de ce que la série, qui a $A\frac{\sin nx}{n}$ pour terme
général, est convergente, l'auteur conclut que la série trigonométrique
générale l'est également. Mais cette conclusion n'est pas permise, car il
est facile de s'assurer que deux séries (du moins lorsque, comme il arrive
ici, les termes n'ont pas tous le même signe) peuvent être l'une
convergente, l'autre divergente, quoique le rapport de deux termes de
même rang diffère aussi peu que l'on veut de l'unité prise positivement
lorsque les termes sont d'un rang très avancé.

On en voit un exemple très simple dans les deux séries, ayant
l'une pour terme général $\frac{(-1)^n}{\sqrt{n}}$, et l'autre
$\frac{(-1)^n}{\sqrt{n}}\bigl(1+\frac{(-1)^n}{\sqrt{n}}\bigr)$. La première
de ces séries est convergente, la seconde au contraire est divergente,
car en la soustraiant de la première on obtient la série divergente:
\[
-1-\tfrac12-\tfrac13-\tfrac14-\tfrac15-\;\text{etc.}
\]
et cependant le rapport se deux termes correspondans, qui est
$1\pm\frac{1}{\sqrt{n}}$,
converge vers l'unité à mesure que $n$ croît.

Je vais maintenant entrer en matière, en commen\c{c}ant par l'examen
des cas les plus simples, auxquels tous les autres peuvent être
ramenés. Désignons par $h$ un nombre positif inférieur ou tout au plus égal
à $\frac\pi2$ et par $f(\beta)$ une fonction de $\beta$ qui reste continue
entre les limites $0$ et~$h$;
j'entends par là une fonction qui a une valeur finie et determinée pour
toute valeur de $\beta$ comprise entre $0$ et $h$, et en outre telle que la
différence $f(\beta+\varepsilon)-f(\beta)$ diminue sans limite lorsque
$\varepsilon$ devient de plus en
plus petit. Supposons encore que la fonction reste toujours positive entre
les limites $0$ et $h$ et qu'elle décroisse constamment depuis $0$ jusquà~$h$,
en sorte que si $p$ et $q$ désignent deux nombres compris entre $0$ et~$h$,
$f(p)-g(q)$ ait toujours un signe opposé à celui de $p-q$. Celà posé
considérons l'intégrale
\[
\tag{1.}
\int_0^h \frac{\sin i\beta}{\sin\beta} f(\beta)\,\partial\beta
\]
dans laquelle $i$ est une quantité positive, et voyons ce que cette intégrale
deviendra à mesure que $i$ croît. Pour cela partageons la en plusieurs
autres prises la première depuis $\beta=0$ jusqu'à $\beta=\frac{\pi}{i}$,
la seconde depuis
$\beta=\frac{\pi}{i}$ jusqu'à $\beta=\frac{2\pi}{i}$
et ainsi de suite, l'avant-dernière ayant
pour limites $(r-1)\frac{\pi}{i}$ et~$\frac{r\pi}{i}$,
et la dernière $\frac{2\pi}{i}$ et~$h$, $\frac{r\pi}{i}$ designant le
plus grand multiple de $\frac{\pi}{i}$ qui soit contenu dans~$h$.
Il est facile de voir
que ces intégrales nouvelles, dont le nombre est $r+1$, sont alternativement
positives et négatives, la fonction placée sous le signe somme étant
évidemment toujours positive entre les limites de la première, négative
entre les limites de la seconde et ainsi de suite. Il n'est pas moins facile
de se convaincre que chacune d'elles est plus petite que la précédente,
abstraction faite du signe. En effet $\nu$ désignant un entier $<r$,
les expressions
\[
\int_{(\nu-1)\frac{\pi}{i}}^{\frac{\nu\pi}{i}}
\frac{\sin i\beta}{\sin\beta} f(\beta)\,\partial\beta
\qquad \text{et} \qquad
\int_{\frac{\nu\pi}{i}}^{(\nu+1)\frac{\pi}{i}}
\frac{\sin i\beta}{\sin\beta} f(\beta)\,\partial\beta
\]
représentent deux intégrales consécutives. Rempla\c{c}ons dans la seconde
$\beta$ par~$\frac{\pi}{i}+\beta$; elle se changera ainsi en celle-ci:
\[
\int_{(\nu-1)\frac{\pi}{i}}^{\frac{\nu\pi}{i}}
\frac{\sin (i\beta+\pi)}{\sin\bigl(\beta+\frac{\pi}{i}\bigr)}
f\Bigl(\beta+\frac{\pi}{i}\Bigr)\,\partial\beta
\]
ou ce qui revient au même:
\[
-\int_{(\nu-1)\frac{\pi}{i}}^{\frac{\nu\pi}{i}}
\frac{\sin i\beta}{\sin\bigl(\beta+\frac{\pi}{i}\bigr)}
f\Bigl(\beta+\frac{\pi}{i}\Bigr)\,\partial\beta.
\]
Les deux intégrales qu'il s'agit de comparer ayant ainsi les mêmes limites,
on voit sans peine que la seconde a une valeur numérique inférieure
à celle de la première. Il suffit pour cela de remarquer qu'il
suit de la supposition que nous avons faite sur la fonction~$f(\beta)$, que
$f\bigl(\frac{\pi}{i}+\beta\bigr)<f(\beta)$
et que d'un autre côté
$\sin\bigl(\frac{\pi}{i}+\beta\bigr)>\sin(\beta)$,
les arcs $\beta$ et~$\frac{\pi}{i}+\beta$
étant l'un et l'autre moindres que~$\frac\pi2$,
car il en résulte l'inégalite
$\frac{f(\beta)}{\sin\beta}>\frac{f(\beta+\frac{\pi}{i})}
{\sin(\beta+\frac{\pi}{i})}$,
qui ayant lieu pour toutes les valeurs de $\beta$
intermédiaires entre les limites $(\nu-1)\frac{\pi}{i}$ et~$\frac{\nu\pi}{i}$,
fait voir que, comme nous
l'avons dit, chaque intégrale est plus grande que celle qui la suit, 
abstraction faite du signe. Cette circonstance a lieu a fortiori,
lorsqu'on compare l'avant-dernière à la dernière, attendu que la différence
des limites $\frac{r\pi}{i}$ et~$h$
de la dernière est inférieure à $\frac{\pi}{i}$ différence commune des limites
de toutes les autres.

Examinons actuellement un peu plus en détail l'intégrale du rang~$\nu$,
qui est
\[
\int_{(\nu-1)\frac{\pi}{i}}^{\frac{\nu\pi}{i}}
\frac{\sin (i\beta)}{\sin\beta}
f(\beta)\,\partial\beta.
\]

Comme la fonction de $\beta$ qui se trouve sous le signe intégral est
le produit des facteurs~$\frac{\sin (i\beta)}{\sin\beta}$, et~$f(\beta)$,
qui sont l'un et l'autre des fonctions
continues de $\beta$ entre les limites de 1'intégration et comme d'un
autre côté le premier de ces facteurs conserve toujours le même signe
entre ces mêmes limites, on conclura en vertu d'un théorème connu, que
l'intégrale que nous considérons est égale à l'intégrale du premier facteur
multipliée par une quantité comprise entre la valeur la plus grande
et la valeur la plus petite de l'autre facteur. Le second facteur décroissant
depuis la première limite jusqu'à la seconde, la quantité dont il s'agit
est comprise entre $f\bigl(\frac{(\nu-1)\pi}{i}\bigr)$
et~$f\bigl(\frac{\nu\pi}{i}\bigr)$.  En la désignant par~$\varrho_\nu$,
notre intégrale sera équivalente à
\[
\varrho_\nu\int_{(\nu-1)\frac{\pi}{i}}^{\frac{\nu\pi}{i}}
\frac{\sin (i\beta)}{\sin\beta}\,\partial\beta.
\]

L'intégrale que renfermée encore cette expression, dépend à la fois de
$\nu$ et de~$i$.
Elle est positive ou négative selon que $v-1$ est pair ou impair;
nous la désignerons désormais par $K_\nu$, abstraction faite du signe.
Nous aurons bientôt besoin de connaître la limite vers laquelle elle
converge, lors
que, $\nu$ restant invariable, $i$ devient de plus en plus grand.
Pour découvrir
cette limite, rempla\c{c}ons $\beta$ par~$\frac{\gamma}{i}$,
$\gamma$ étant une nouvelle variable. Nous
aurons ainsi
\[
\int_{(\nu-1)\pi}^{\nu\pi}
\frac{\sin\gamma}{i\sin(\frac{\gamma}{i})}\,\partial\gamma.
\]
Sous cette forme, il est évident qu'elle converge vers la limite
\[
\int_{(\nu-1)\pi}^{\nu\pi}
\frac{\sin\gamma}{\gamma}\,\partial\gamma,
\]
que pour abréger nous désignerons par~$k_\nu$, abstraction faite du signe.

On sait que l'intégrale
$\int_0^\infty\frac{\sin\gamma}{\gamma}\,\partial\gamma$
a une valeur finie et égale à~$\frac\pi2$.
Cette intégrale peut être partagée en une infinité d'autres, prises la première
depuis $\gamma=0$ jusqu'à~$\gamma=\pi$, la seconde depuis $\gamma=\pi$ jusqu'à
$\gamma=2\pi$, et ainsi de suite. Ces nouvelles intégrales sont alternativement
positives et négatives, chacune d'elles a une valeur numérique inférieure
à celle de la précédente, et celle du rang~$\nu$ est~$k_\nu$,
abstraction faite du
signe. La proposition qu'on vient de citer, revient donc à dire que la
suite infinie
\[
\tag{2.}
k_1-k_2+k_3-k_4+k_5-\;\text{etc.}
\]
est convergente et a une somme égale à~$\frac\pi2$.

Les termes de cette suite allant toujours en décroissant, il suit
d'une proposition connue que la somme de $n$ premiers termes est supérieure
ou inférieure à~$\frac\pi2$, selon que $n$ est impair ou pair et que cette
somme qu'on peut désigner par~$S_n$, diffère de $\frac\pi2$ d'une quantité 
moindre que le terme suivant~$k_{n+1}$.

Reprenons actuellement l'intégrale (1.) et cherchons à déterminer
la limite vers laquelle elle converge lorsque $i$ croît indéfiniment. En
faisant ainsi croître le nombre~$i$, les intégrales dans lesquelles nous avons
décomposé l'intégrale~(1.), changeront sans cesse de valeur en même
temps que leur nombre augmentera; il s'agit de connaître le résultat de
ce double changement lorsqu'il continue indéfiniment. Pour cela, prenons
un nombre entier~$m$ (qu'il soit supposé pair pour plus de simplicité)
et supposons que le nombre $m$ reste invariable pendant que $i$ croît.
Le nombre~$r$, qui croît sans cesse avec~$i$, finira bientôt par surpasser le
nombre invariable~$m$, quelque grand qu'on l'ait choisi.

Cela posé, partageons en deux groupes les intégrales dont la somme
est équivalente à l'intégrale~(1.). Le premier groupe comprendra les~$m$
premières de ces intégrales, et le second sera composé de toutes les suivantes.
On aura pour la somme du premier groupe:
\[
\tag{3.}
K_1\varrho_1-K_2\varrho_2+K_3\varrho_3-K_4\varrho_4+\cdots-K_m\varrho_m
\]
et le second, dont le nombre des termes croît sans cesse avec~$i$, a pour
premiers termes:
\[
\tag{4.}
K_{m+1}\varrho_{m+1}-K_{m+2}\varrho_{m+2}+\cdots.
\]

Considérons séparément ces deux groupes. Le nombre $i$ croissant
indéfiniment la somme~(3.) convergera vers une limite qu'il est facile de
déterminer. En effet, les quantités $\varrho_1,\varrho_2,\ldots\varrho_m$
qui sont comprises la
première entre $f(0)$ et~$f(\frac{\pi}{i})$, la seconde entre
$f(\frac{\pi}{i})$ et~$f(\frac{2\pi}{i})$, et la dernière entre
$f(\frac{(m-1)\pi}{i})$ et $f(\frac{m\pi}{i})$
convergent chacune vers la limite~$f(0)$,
lorsque, $m$ restant invariable, $i$ croît sans cesse. D'un autre côté nous
avons vu que les quantités $K_1,K_2,\ldots K_m$ convergent dans les mêmes
circonstances respectivement vers les limites $k_1,k_2,\ldots k_m$. Donc
la somme (3.) converge vers la limite $f(0)(k_1-k_2+k_3-\text{etc.}\ldots
-k_m)=S_mf(0)$,
ce qui veut dire que la différence entre la somme (3.) et
$S_mf(0)$ finira toujours, abstraction faite du signe, par être constamment
inférieure à~$\omega$, $\omega$ désignant une quantité positive aussi
petite que l'on veut.

Considérons pareillement la somme (4.), dont le nombre des termes
augmente sans cesse. Ses termes étant alternativement positifs et
négatifs, et chacun d'eux ayant une valeur numérique inférieure à celle
du terme précédent, comme nous l'avons vu plus haut, en considérant
les intégrales que ces termes représentent, il suit d'un principe
connu\footnote{Le principe sur lequel nous nous appuyons peut être
énoncé de cette manière. Les lettres $A,A',A'',\ldots$  désignant des
quantités positives en nombre quelconque et telles que
\[
A>A'>A''>\text{etc., la quantité\ } A-A'+A''-A'''+\text{etc.}
\]
est positive et inférieure à~$A$. Celà résulte immédiatement de que
la quantité précedente peut être
mise sous l'une et l'autre de ces deux formes:
\[
(A-A') + (A''-A''')+\text{etc.,}
\]
\[
A-(A'- A'') - (A'''- A^\text{IV})-\text{etc.}
\]},
que cette somme, quelque soit le nombre de ses termes, est positive comme
son premier terme~$K_{m+1}\varrho_{m+1}$, et a une valeur inférieure
à celle de ce
terme. Or, ce premier terme convergeant vers la limite~$k_{m+1}f(0)$, il
s'ensuit que la somme~(4.) finira toujours par être inférieure à
$k_{m+1}f(0)$
augmenté d'une quantité positive $\omega'$ aussi petite que l'on veut.
En combinant ce résultat avec celui que nous avons obtenu sur la somme~(3.),
il n'y a qu'un instant, on verra que l'intégrale~(1.) qui est la somme des
expressions (3.) et~(4.) finira toujours par différer de $f(0)S_m$
d'une quantité moindre, abstraction faite du signe, que
$\omega+\omega'+f(0)k_{m+1}$, $\omega$ et $\omega'$
étant deux nombres d'une petitesse arbitraire. D'un autre côté $S_m$
diffère de $\frac\pi2$ d'une quantité numériquement inférieure à~$k_{m+1}$; 
donc l'intégrale finira toujours par différer de
$\frac\pi2f(0)$ d'une quantité moindre que
$\omega+\omega'+2f(0)k_{m+1}$, abstraction faite du signe.

Comme $m$ peut être choisi tellement grand que~$k_{m+1}$, sait moindre
que toute grandeur donnée, il s'ensuit que l'intégrale~(1.) finira toujours,
lorsque $i$ croît sans limite, par différer constamment de $\frac\pi2f(0)$
d'une quantité moindre, abstraction faite du signe, qu'un nombre aussi
petit que l'on veut. Il est ainsi prouvé, que l'intégrale~(1.) converge vers
la limite $\frac\pi2f(0)$ pour des valeurs croissantes de~$i$.

Supposons maintenant que la fonction $f(\beta)$, au lieu d'être toujours
décroissante depuis $0$ jusqu'à~$h$, soit constante et égale à l'unité.
On pourra
dans ce cas déterminer la limite vers laquelle converge l'intégrale~(1.)
par les mêmes considérations que nous venons d'employer; c'est ce qu'on
voit tout de suite, en se rappellant que la démonstration précédente est
fondée sur ce que les intégrales dans lesquelles nous avons decomposé
l'intégrale~(1.), forment une suite décroissante. Or, ce décroissement tient
à deux choses, au décroissement du facteur~$f(\beta)$ et à l'accroissement du
diviseur~$\sin\beta$. Si $f(\beta)$ devient un nombre constant,
l'accroissement de $\sin\beta$ suffira toujours pour rendre chaque
intégrale de la série plus petite
que la précédente. On trouvera ainsi, en supposant toujours $h$ positive
et tout au plus égale à~$\frac\pi2$, que l'intégrale
$\int_0^h \frac{\sin i\beta}{\sin\beta}\,\partial\beta$ converge vers la
limite~$\frac\pi2$.   Il suit de là que l'intégrale
$\int_0^h \frac{c\sin i\beta}{\sin\beta}\,\partial\beta$, dans laquelle $c$ est
une constante positive ou négative, converge vers la limite~$c\frac\pi2$.

Nous avons supposé que la fonction $f(\beta)$ était décroissante et positive
entre les limites $0$ et~$h$. La première circonstance ayant toujours
lieu, c'est-à-dire la fonction étant telle que $f(p)-f(q)$ ait un signe 
contraire à celui de $p-q$ pour des valeurs $p$ et $q$ comprises entre $0$
et~$h$,
supposons que $f(\beta)$ ne soit pas toujours positive. On prendra
alors une constante positive $c$ assez grande pour que
$c+f(\beta)$ conserve toujours un signe
positif depuis $\beta = 0$ jusqu'à~$\beta=h$.
L'intégrale $\int_0^h f(\beta)\frac{\sin i\beta}{\sin\beta}\,\partial\beta$
étant égale à la différence de celles-ci:
$\int_0^h [c+f(\beta)]\frac{\sin i\beta}{\sin\beta}\,\partial\beta$ et
$\int_0^h c\frac{\sin i\beta}{\sin\beta}\,\partial\beta$,
sa limite sera la différence des limites vers lesquelles convergent ces 
dernières. Or ces dernières rentrent dans les cas précédemment examinés
($c+f(\beta)$ étant une fonction décroissante et positive) et convergent vers
les limites $[c+f(0)]\frac\pi2$ et~$c\frac\pi2$, d'où il suit que
la première converge vers la limite~$\frac\pi2f(0)$.

Considérons actuellement une fonction $f(\beta)$ croissante depuis $0$
jusqu'à~$h$. Dans ce cas $-f(\beta)$ sera une fonction décroissante.
L'intégrale $\int_0^h -f(\beta)\frac{\sin i\beta}{\sin\beta}\,\partial\beta$
convergera donc vers la limite~$-\frac\pi2f(0)$, et par conséquent l'intégrale
$\int_0^h f(\beta)\frac{\sin i\beta}{\sin\beta}\,\partial\beta$
vers la limite~$\frac\pi2f(0)$.

En réunissant ces résultats, on aura cet énoncé:
\[
\tag{5.}
\left\{\parbox{0.85\textwidth}{
»Quelle que soit la fonction $f(\beta)$, pourvu qu'elle reste continue
entre les limites $0$ et~$h$ ($h$ étant positive et tout au plus égale
à~$\frac\pi2$),
et qu'elle croisse ou qu'elle décroisse depuis la première de ces limites
jusqu'à la seconde, l'1ntégrale
$\int_0^h f(\beta)\frac{\sin i\beta}{\sin\beta}\,\partial\beta$
finira par différer constamment de
$\frac\pi2$ d'une quantité moindre que tout nombre
assignable, lorsqu'on y fait croître $i$ au delà de toute limite positive.«}
\right.
\]

Désignons par $g$ un nombre positif différent de zéro et inférieur
à~$h$, et supposons que la fonction reste continue et croisse ou décroisse
depuis $g$ jusquà~$h$. L'intégrale
$\int_0^h f(\beta)\frac{\sin i\beta}{\sin\beta}\,\partial\beta$
convergera alors vers
une limite qu'il est facile de découvrir. On pourroit y parvenir par des
considérations analogues à celles que nous avons appliquées à l'intégrale~(1.);
mais il est plus simple de ramener ce nouveau cas à ceux que nous
avons considérés dans ce qui précède. La fonction n'étant donnée que
depuis $g$ jusqu'à $h$ reste entièrement arbitraire pour les valeurs de $\beta$
comprises entre $0$ et~$g$. Supposons que l'on entende par~$f(\beta)$, pour les
valeurs de $\beta$ comprises entre $0$ et $g$ une fonction continue
et croissante
ou décroissante depuis $0$ jusqu'à~$g$, selon que $f(\beta)$
croît ou décroît depuis
$g$ jusqu'à~$h$; supposons encore que $f(g-\varepsilon)$ diffère infiniment
peu de~$f(g+\varepsilon)$, si $\varepsilon$
décroît sans limite; ayant satisfait d'une manière quelconque à
ces conditions, ce qu'on peut toujours faire d'une infinité de manières,
la fonction $f(\beta)$ remplira depuis $0$ jusqu'à $h$ les conditions
exprimées dans l'énoncé~(5.). Les intégrales
\[
\int_0^g f(\beta)\frac{\sin i\beta}{\sin\beta}\,\partial\beta
\qquad \text{et} \qquad
\int_0^h f(\beta)\frac{\sin i\beta}{\sin\beta}\,\partial\beta
\]
convergeront donc l'une et l'autre vers la limite~$\frac\pi2f(0)$. D'où l'on 
conclut que l'intégrale
$\int_g^h f(\beta)\frac{\sin i\beta}{\sin\beta}\,\partial\beta$
qui est la différence des précédentes,
a zéro pour limite.

Ce nouveau résultat peut être réuni en un seul énoncé avec celui
que nous avons obtenu plus haut. On aura ainsi:
\[
\tag{6.}
\left\{\parbox{0.85\textwidth}{
»La lettre $h$ désignant une quantité positive tout au plus égale à~$\frac\pi2$,
et $g$ une quantité également positive et en outre inférieure à~$h$,
l'intégrale
\[
\int_g^h f(\beta)\frac{\sin i\beta}{\sin\beta}\,\partial\beta
\]
dans laquelle la fonction $f(\beta)$ est continue entre les limites de
l'intégration et a une marche toujours croissante ou toujours décroissante
depuis $\beta=g$ jusqu'à~$\beta=h$, convergera vers une certaine
limite, lorsque le nombre $i$ devient de plus en plus grand. Cette
limite est égale à zéro, le seul cas excepté où $g$ a une valeur nulle,
dans ce cas elle a la valeur~$\frac\pi2f(0)$.«}
\right.
\]

Il est évident que ce résultat ne serait que légèrement modifié, si la
fonction $f(\beta)$ présentait une solution de continuité pour~$\beta=g$,
et $\beta=h$,
c'est-à-dire si $f(g)$ était différent de $f(g+\varepsilon$) et $f(h)$
de~$f(h-\varepsilon)$, $\varepsilon$ désignant
une quantité infiniment petite et positive, pourvu qu'alors les valeurs
$f(g)$ et $f(h)$ ne fûssent pas infinies. Il faudrait seulement dans ce
cas remplacer $f(0)$ par $f(\varepsilon)$ dans l'énoncé précédent,
ce qu'on peut faire
encore même quand il n'y a pas de solution de continuité, attendu qu'alors
$f(\varepsilon)$ est égale à~$f(0)$.

Nous sommes maintenant en état de prouver la convergence des
séries périodiques qui expriment des fonctions arbitraires entre des limites
données. La marche que nous allons suivre nous conduira à établir
la convergence de ces séries et à déterminer en même temps leurs valeurs.
Soit $\phi(x)$ une fonction de~$x$, ayant une valeur finie et déterminée
pour chaque valeur de $x$ comprise entre $-\pi$ et~$\pi$, et supposons qu'il
s'agisse de développer cette fonction dans une série de sinus et de cosinus
d'arcs multiples de~$x$. La série qui résout cette question, est, comme
l'on sait:
\[
\arraycolsep=1.5pt
\tag{7.}
\frac{1}{2\pi}\int\phi(\alpha)\,\partial\alpha+\frac1\pi
\left\{\begin{array}{l}
\displaystyle \cos x\int\phi(\alpha)\cos\alpha\,\partial\alpha+
\cos2x\int\phi(\alpha)\cos2\alpha\,\partial\alpha\cdots\\
\displaystyle \sin x\int\phi(\alpha)\sin\alpha\,\partial\alpha+
\sin2x\int\phi(\alpha)\sin2\alpha\,\partial\alpha\cdots
\end{array}\right\}.
\]
Les intégrales qui déterminent les coëfficiens constans, étant prises depuis
$\alpha=-\pi$ jusqu'à~$\alpha=\pi$, et $x$ désignant une quantité quelconque
comprise entre $-\pi$ et~$\pi$ (\Title{Théorie de la Chaleur},
No.~232. et suiv.).

Considérons les $2n+1$ premiers termes de cette série ($n$ étant un
nombre entier) et voyons vers quelle limite converge la somme de ces
termes, lorsque $n$ devient de plus en plus grand. Cette somme peut être
mise sous la forme suivante:
\[
\frac1\pi\int_{-\pi}^{+\pi}\phi(\alpha)\,\partial\alpha
\bigl[\tfrac12+\cos(\alpha-x)+\cos2(\alpha-x)+\cdots+\cos n(\alpha-x)\bigr],
\]
ou en sommant la suite de cosinus,
\[
\tag{8.} \frac1\pi\int_{-\pi}^{+\pi}\phi(\alpha)
\frac{\sin(n+\frac12)(\alpha-x)}{2\sin\frac12(\alpha-x)}\,\partial\alpha.
\]

Tout se réduit maintenant à déterminer la limite dont cette intégrale
approche sans cesse, lorsque $n$ croît indéfiniment. Pour cela nous
la partagerons en deux autres prises l'une depuis $-\pi$ jusqu'à~$\pi$, l'autre
depuis $x$ jusqu'à~$\pi$. Si l'on remplace dans la première $\alpha$
par~$x-2\beta$,
et dans la seconde $\alpha$ par~$x+2\beta$, $\beta$ étant une nouvelle
variable, ces deux
intégrales se changeront en celles-ci, abstraction faite du 
facteur~$\frac1\pi$: 
\[
\tag{9.,10.}
\int_0^{\frac12(\pi+x)}
\frac{\sin(2n+1)\beta}{\sin\beta}\phi(x-2\beta)\,\partial\beta
\quad\text{et\ }
\int_0^{\frac12(\pi-x)}
\frac{\sin(2n+1)\beta}{\sin\beta}\phi(x+2\beta)\,\partial\beta.
\]

Considérons la seconde de ces deux intégrales. La quantité $x$
étant inférieure ou tout au plus égale à~$\pi$, abstraction faite du signe,
$\frac12(\pi-x)$
ne pourra tomber hors des limites $0$ et~$\pi$. Si~${\frac12(\pi-x)=0}$,
ce qui a lieu lorsque~${x=\pi}$, l'intégrale est nulle quelque que soit~$n$;
dans tous les autres cas elle convergera pour des valeurs croissantes de $n$
vers une limite que nous allons déterminer. Supposons d'abord $\frac12(\pi-x)$
inférieure ou tout au plus égale à~$\frac\pi2$, et remarquons que la fonction
$\phi(x+2\beta)$ peut présenter plusieurs solutions de continuité depuis
$\beta=0$ jusqu'à~${\beta=\frac12(\pi-x)}$,
et qu'elle peut aussi avoir plusieurs maxima et
minima dans ce même intervalle. Désignons par $l,l',l'',\ldots l^{(\nu)}$
rangées selon l'ordre de leur grandeur, les différentes valeurs de~$\beta$,
qui présentent l'une ou l'autre de ces circonstances, et décomposons notre
intégrale en plusieurs autres prises respectivement entre les limites $0$ 
et~$l$, $l'$ et~$l''$, $\ldots l^{(\nu)}$ et~$\frac12(\pi-x)$.
Toutes ces intégrales se trouveront dans
le cas de l'énoncé~(6.). Elles convergeront donc toutes vers la limite zéro
à mesure que $n$ croît, à l'exception de la première qui converge vers la
limite~$\frac\pi2\phi(x+\varepsilon)$, $\varepsilon$
étant un nombre infiniment petit et positif. Si
$\frac12(\pi-x)$ était superieure à~$\frac\pi2$,
ce qui arrivera lorsque $x$ a une valeur
négative, on partagerait l'intégrale en deux autres, l'une prise depuis
$\beta=0$ jusqu'à~$\beta=\frac\pi2$, l'autre depuis $\frac\pi2$
jusqu'à~$\beta=\frac12(\pi-x)$. La première
de ces nouvelles intégrales se trouvera dans le même cas que celle
que nous venons de considérer, elle convergera donc vers la
limite~$\frac\pi2\phi(x+\varepsilon)$.
Quant à la seconde, on peut la changer en celle-ci, en y
rempla\c{c}ant $\beta$ par~$\pi-\gamma$, $\gamma$
étant une nouvelle variable:
\[
\int_{\frac12(\pi+x)}^{\frac\pi2}\phi(x+2\pi-2\gamma)
\frac{\sin(2n+1)(\pi-\gamma)}{\sin(\pi-\gamma)}\,\partial\gamma,
\]
ou ce qui revient au même, $n$ étant un entier:
\[
\int_{\frac12(\pi+x)}^{\frac\pi2}\phi(x+2\pi-2\gamma)
\frac{\sin(2n+1)\gamma}{\sin\gamma}\,\partial\gamma.
\]

Elle a ainsi une forme analogue à celle de la précédente; en la décomposant
comme elle en plusieurs autres, on verra qu'elle converge
vers la limite zéro, le seul cas excepté, où $\frac12(\pi+x)$
a une valeur nulle,
c'est-à-dire lorsque~$x=-\pi$; dans ce cas elle approche continuellement
de la limite~$\phi(\pi-\varepsilon)$, $\varepsilon$
ayant toujours la même signification. En résumant
tout ce qui précède, on trouvera que la seconde des intégrales (10.)
est nulle lorsque~$x=\pi$, qu'elle converge vers la limite
$\frac\pi2\bigl[\phi(\pi-\varepsilon)+\phi(-\pi+\varepsilon)\bigr]$
lorsque~$x=-\pi$, et que dans tous les autres
cas elle approche continuellement de la limite~$\frac\pi2\phi(x+\varepsilon)$.
La première des intégrales (9.) est entièrement analogue
à la seconde; en y appliquant
des considérations semblables, on trouvera qu'elle est nulle
lorsque~$x=-\pi$, qu'elle converge vers la limite
$\frac\pi2\bigl[\phi(\pi-\varepsilon)+\phi(-\pi+\varepsilon)\bigr]$
lorsque $x=\pi$ et
que dans tous les autres cas elle a pour limite~$\frac\pi2\phi(x-\varepsilon)$.
Connoissant ainsi
les limites de chacune des intégrales~(9.,10.),
il est facile de trouver la limite
dont l'intégrale (8.) approche sans cesse, lorsque $n$ devient de plus en plus
grand; il suffit pour celà de se rappeler que cette intégrale est égale à
la somme des intégrales (9.,10.) divisée par~$\pi$. Or, l'intégrale (8.) étant
équivalente à la somme des $2n+1$ premiers termes de la série~(7.), il
est prouvé que cette série est convergente et l'on trouve au moyen des
résultats précédens qu'elle est égale à
$\frac12\bigl[\phi(x+\varepsilon)-\phi(x-\varepsilon)\bigr]$
pour toute
valeur de $x$ comprise entre $-\pi$ et~$\pi$, et que pour chacune des valeurs
extrêmes $\pi$ et~$-\pi$, elle est égale
à~$\frac12\bigl[\phi(\pi-\varepsilon)+\phi(-\pi+\varepsilon)\bigr]$.

L'exposé précédent embrasse tous les cas; il se simplifie lorsque
la valeur de $x$ qu'on considère n'est pas une de celles qui présentent
une solution de continuité. En effet les quantités
$\phi(x+\varepsilon)$ et $\phi(x-\varepsilon)$
étant alors l'une et l'autre équivalentes à~$\phi(x)$, on voit que Ia série a
pour valeur~$\phi(x)$.

Les considérations précédentes prouvent d'une manière rigoureuse
que, si la fonction $\phi(x)$, dont toutes les valeurs sont supposées finies et
déterminées, ne présente qu'un nombre fini de solutions de continuité entre
les limites $-\pi$ et~$\pi$, et si en outre elle n'a qu'un nombre déterminé
de maxima et de minima entre ces mêmes limites, la série (7.),
dont les coëfficiens sont des intégrales définies dépendantes de la fonction
$\phi(x)$ est convergente et a une valeur généralement exprimée par
$\frac12\bigl[\phi(x+\varepsilon)+\phi(x-\varepsilon)\bigr]$,
où $\varepsilon$ désigne un nombre infiniment petit. Il nous
resterait à considérer les cas où les suppositions que nous avons faites
sur le nombre des solutions de continuité et sur celui des valeurs maxima
et minima cessent d'avoir lieu. Ces cas singuliers peuvent être ramenés
à ceux que nous venons de considérer. ll faut seulement pour que la
série (8.) présente un sens lorsque les solutions de continuité sont en
nombre infini, que la fonction $\phi(x)$ remplisse la condition suivante.

Il est nécessaire qu'alors la fonction $\phi(x)$ soit telle que, si l'on 
désigne par $a$ et $b$ deux quantités quelconques comprises entre
$-\pi$ et~$\pi$,
on puisse toujours placer entre $a$ et~$b$ d'autres quantités
$r$ et $s$ assez rapprochées
pour que la fonction reste continue dans l'intervalle de $r$ à~$s$.
On sentira facilement la nécessité de cette restriction en considérant que
les différens termes de la série sont des intégrales définies et en remontant
à la notion fondamentale des intégrales. On verra alors que l'intégrale
d'une fonction ne signifie quelque chose qu'autant que la fonction
satisfait à la condition précédemment énoncée. On aurait un exemple
d'une fonction qui ne remplit pas cette condition, si l'on supposait $\phi(x)$
égale à une constante déterminée $c$ lorsque la variable $x$ obtient une
valeur rationnelle, et égale à une autre constante $d$, lorsque cette variable
est irrationnelle. La fonction ainsi définie a des valeurs finies et
déterminées pour toute valeur de~$x$, et cependant on ne saurait la substituer
dans la série, attendu que les différentes intégrales qui entrent dans cette
série, perdroient toute signification dans ce cas. La restriction que je
viens de préciser, et celle de ne pas devenir infinie, sont les seules 
auxquelle la fonction $\phi(x)$ soit sujette
et tous les cas qu'elles n'excluent pas
peuvent être ramenés à ceux que nous avons considérés dans ce qui
précède. Mais la chose, pour être faite avec toute la clarté qu'on peut
désirer exige quelques détails liés aux principes fondamentaux de l'analyse
infinitésimale et qui seront exposés dans une autre note, dans laquelle
je m'occuperai aussi de quelques autres propriétés assez remarquables
de la série~(7.).

$\qquad$ Berlin, Janvier 1829.

\end{document}